\documentclass[a4paper,fleqn]{article}  
\usepackage{tikz}
\usepackage{url}
\usepackage{float}  
\usepackage{placeins} 

\usetikzlibrary{arrows.meta, positioning, shapes.misc, shadows.blur}
\usepackage{graphicx}
\usepackage{caption}
\usepackage{subcaption}

\usepackage{pgfplots}         
\usepackage{pgfplotstable}    
\usepackage{amsmath, amssymb} 

\usepackage{siunitx}          
\usepackage{xcolor}           

\usepackage{booktabs}         

\pgfplotsset{compat=1.18}     

\usepackage[utf8]{inputenc}
\usepackage[T1]{fontenc}
\usepackage{lmodern}
\usepackage{geometry}
\usepackage{setspace}
\usepackage{authblk} 
\usepackage{amsmath, amssymb}

\usepackage[authoryear,longnamesfirst]{natbib}

\title{\vspace{-2cm}\textbf{Is simplicity still possible for a more accurate approximation to the perimeter of the ellipse? or, Using the exponential function to further improve the second Ramanujan's approximation}}
\author[1]{Salvador E. Ayala-Raggi\thanks{Corresponding author: \texttt{salvador.raggi@correo.buap.mx}}}
\author[1]{Manuel Rendón-Marín\thanks{\texttt{manuel.rendon@correo.buap.mx}}}
\affil[1]{Facultad de Ciencias de la Electrónica, Benemérita Universidad Autónoma de Puebla (BUAP), Puebla, México}
\date{\vspace{0.3cm}October 2025\\[1em]\small Preprint submitted to arXiv}

\newcommand{\keywords}[1]{%
  \vspace{0.5em}
  \noindent\textbf{Keywords: }#1
  \vspace{1em}
}

\begin{document}

\vspace{-1cm}

\maketitle

\begin{abstract}
The perimeter of an ellipse has no exact closed-form expression in terms of elementary functions, and numerous approximations have been proposed since the eighteenth century. Classical formulas by Fagnano, Euler, and Ramanujan, as well as modern refinements such as Cantrell and Koshy methods, aim to reduce the approximation error while maintaining computational simplicity. In this paper, we introduce a new closed-form expression that enhances Ramanujan’s second formula by dividing it by 1 minus a binomial of two exponential terms resulting in a very stable approximation in a range of $b/a$ between $1$ and $0$. The resulting approximation remains compact, requiring only four constants, and achieving a remarkable tradeoff between simplicity and accuracy. Across the full eccentricity range of $e \in [0,1]$, our method attains a maximum relative error of approximately $0.57$~ppm with respect to the exact perimeter computed via elliptic integral. Our formula is quasi-exact at the extremes, for the circle ($b/a=1$) and for the degenerate flat ellipse ($b/a=0$). Compared with Ramanujan2-Cantrell approximation, the proposed method reduces the maximum relative error by a factor of $25$ while preserving a short and elegant expression. This makes it one of the simplest yet most accurate closed-form and single-line approximations to the ellipse perimeter currently available in the literature.
\end{abstract}


\section{Introduction}

The exact computation of the perimeter of an ellipse has been a classical problem in mathematics since antiquity. Unlike the circumference of a circle, which admits a closed-form expression $P=2\pi r$, the ellipse with semimajor axis $a$ and semiminor axis $b$ $(a\geq b)$ does not have a simple algebraic expression for its perimeter. Instead, the exact formula involves the complete elliptic integral of the second kind, $E(e)$, where $e$ is the eccentricity, i.e.,
\begin{equation}
P_{\mathrm{exact}} = 4a\,E\!\left(e\right), \quad \text{where} \quad E(e) = \int_0^{\pi/2} \sqrt{1 - e^2 \sin^2\theta}\, d\theta \quad \text{and} \quad e=\sqrt{1-\frac{b^2}{a^2}}.
\label{eq:exact}
\end{equation}
Because elliptic integrals resisted closed-form evaluation for centuries, mathematicians and scientists proposed numerous approximations, balancing simplicity with accuracy. We summarize seven formulas of historical and practical importance.

\subsection*{Fagnano's Formula}
One of the earliest attempts is attributed to Giulio Carlo de’ Fagnano (18th century) \cite{Fagnano1735}, \cite{Sykora2005} . He proposed the very simple approximation
\begin{equation}
P_{\text{Fagnano}} = \pi (a+b).
\end{equation}
While exact for the circle ($a=b$), this formula becomes increasingly inaccurate for ellipses with large eccentricity.

\subsection*{Euler's Formula}
Leonhard Euler improved on this by considering the quadratic mean of the semiaxes \cite{Euler1740}, \cite{Sykora2005}, \cite{NemesSykora2007}, giving
\begin{equation}
P_{\text{Euler}} = 2\pi \sqrt{\tfrac{a^2+b^2}{2}}.
\end{equation}
This expression is significantly more accurate than Fagnano's but still yields large errors for highly elongated ellipses.

\subsection*{Euler-Ivory's series expansion}

A classical expansion in infinite series of the elliptic integral was proposed by Euler in 1773 in terms of the eccentricity $e$. To accelerate convergence, it was rewritten in terms of $h$ by Ivory in 1825 \cite{Michon2007}, although it is also known as Gauss-Kummer series of $h$:
\begin{equation}
P_{Euler-Ivory} = \pi (a + b) \sum_{n=0}^{\infty} 
\left[ 
\binom{\tfrac{1}{2}}{n}
\right]^{2} 
h^{n}
= 
\pi (a + b)
\left[
1 + \frac{h}{4} + \frac{h^{2}}{64} + \frac{h^{3}}{256} + \cdots
\right], \qquad \text{with} \qquad h = \left(\frac{a-b}{a+b}\right)^2.
\label{eq:Ivory}
\end{equation}

\subsection*{Ramanujan I}
In the early 20th century, Srinivasa Ramanujan proposed remarkably accurate approximations \cite{Ramanujan1914a}, \cite{Sykora2005}, \cite{NemesSykora2007}. His first formula (R1) is
\begin{equation}
P_{\text{R1}} = \pi \left[ 3(a+b) - \sqrt{(3a+b)(a+3b)} \right].
\end{equation}
This compact formula remains widely cited due to its balance of simplicity and improved accuracy.

\subsection*{Ramanujan II}
Ramanujan also proposed, in terms of $h$, a second even more precise expression (R2) \cite{Ramanujan1914a}, \cite{Sykora2005}, \cite{NemesSykora2007}:
\begin{equation}
P_{\text{R2}} = \pi (a+b) \left( 1 + \frac{3h}{10 + \sqrt{4 - 3h}} \right), \qquad \text{with} \qquad h = \left(\frac{a-b}{a+b}\right)^2.
\label{eq:R2}
\end{equation}
This approximation achieves errors close to $4. \times 10^{-4}$ or less across a wide range of axis ratios.

\subsection*{Cantrell’s Formula}
Building upon Ramanujan II, Cantrell (2004) \cite{Cantrell2004}, \cite{Sykora2005}, \cite{NemesSykora2007} introduced an additional corrective term to further reduce the error:
\begin{equation}
P_{\text{Cantrell}} = \pi (a+b) \left( 1 + \frac{3h}{10 + \sqrt{4 - 3h}} + c h^{12} \right), \quad \text{where} \quad c = \frac{4}{\pi} - \frac{14}{11}.
\end{equation}
This modification reduces the maximum relative error to about $1.4 \times 10^{-5}$, representing one of the best-known simple closed-form approximations.

\subsection*{Koshy’s Formulas for quarter-perimeter}
More recently, Koshy (2024) proposed two formulas for the quarter-perimeter \cite{Koshy2024}, denoted $Q(a,b)$. The first formula, more accurate than the second one, is of the form

\begin{equation}
Q(a,b) \;\approx\; Q(a,b;p) \;=\; \bigl(a^p + b^p \bigr)^{1/p}, \quad \text{where} \quad p = p(a,b;k)=\frac{\ln(2)}{\ln(\pi/2)} + 1 - (b/a)^k, 
\end{equation}
and $k \;=\; 0.03214 \;-\; 0.0734\,(b/a) \;+\; 0.0863\,(b/a)^2 \;-\; 0.0681\,(b/a)^3 \;+\; 0.02306\,(b/a)^4 \;$. Finally, the perimeter can be calculated as $P_{Koshy1} \approx 4 \, Q(a,b)$. This approach yields a maximum error close to 8.7 ppm.
The second formula proposed by Koshy is:
\begin{equation}
Q(a,b) \approx Q(a,b;2,k) = 
\sqrt{a^2 + b^2} 
+ \left( \frac{\pi}{2} - \sqrt{2} \right)
\left( \frac{GM}{AM} \right)^k 
\sqrt{a b},
\end{equation}
where $k = 2.6071 + 1.2243\left(\frac{b}{a}\right) - 1.2673\left(\frac{b}{a}\right)^2 + 0.45566\left(\frac{b}{a}\right)^3$, $AM = \frac{a + b}{2}$, and $GM = \sqrt{a b}$.

\subsection*{Recent approaches}
There are other recent works, like the one proposed by Moscato in \cite{Moscato2024}, giving a maximum absolute relative error of 2.7 ppm. In \cite{Sykora2007Advances}, Sykora describes an interesting Ahmadi's approximation that achieves a maximum relative error < 1 ppm. However, these formulas are not single-line, they are complicated, and particularly in the Ahmadi's case it includes a total of eight exponentiations.  In our work, the focus has been on a much more compact, closed-form and easy-to-calculate expression.


\section{Materials and Methods}


\subsection{One–exponential correction to R2 ($R2/1exp$)}
We started by studying the relative error curve of the classical second Ramanujan formula $P_{\text{R2}}$ given in Eq. (\ref{eq:R2}). From now on, for all formulas the relative error will be calculated as follows:
\begin{equation}
e_{r}(P_{test}) =\frac{P_{test}-P_{exact}}{P_{exact}}.
\label{eq:error}
\end{equation}
Using Eq. (\ref{eq:error}), the relative error of $P_{R2}$  is $e_{r}(P_{R2})$. We have seen that $e_{r}(P_{R2})$ expressed as a function of $h$, is smooth, negative, and increases in magnitude as $h\to 1$ (high eccentricity), resembling a decaying exponential in $1-h$ of the form $\Delta e = -Ae^{-B(1-h)}$. Therefore, if we start from the hypothesis that there must exist $\Delta e$ such that $e_{r} - \Delta e = 0$, then 

\begin{equation}
\hat{P}_{\mathrm{1}} \;=\; \frac{P_{\mathrm{R2}}}{1 - A\,e^{-B(1-h)}}.
\label{eq:corr1}
\end{equation}
The parameters $(A,B)$ are fitted once (and reused universally) by minimizing a uniform error criterion over a mesh. 
Specifically, we minimize the maximum absolute relative error
\begin{equation}
\Phi(A,B)\;=\;\max_{(a,b)\in\mathcal{D}_1}
\frac{\big|\hat{P}_{\mathrm{1}}-P_{\mathrm{exact}}\big|}
     {P_{\mathrm{exact}}},
\qquad
\mathcal{D}_1=\{b=1,\; a\in[1,100]\}.
\end{equation}
In this way we obtain:
\[
A \;=\; 3.62077\times 10^{-4}, 
\qquad 
B \;=\; 10.826,
\]
which yields a maximum relative error of about $2.14$\ ppm on $\mathcal{D}_1$ (see Figure \ref{fig:rel_error_methods}(a)). Equation \ref{eq:corr1} can be used reliably for $a \leq 100$. As an interesting fact, in this interval are all the elliptical orbits of all known comets.

\subsection{Two–exponential correction ($R2/2exp$)}
To further suppress the residual error for very low $b/a$ ratio, i.e. $a > 100$, while preserving compactness, we added a second exponential term of the same shape with parameters $C$ and $D$, which are fitted with the same method as for $A$ and $B$ used in $\hat{P}_{\mathrm{1}}$. The new formula is:
\begin{equation}
\hat{P}_{\mathrm{2}}= \frac{P_{\mathrm{R2}}}{1 - \Big[A\,e^{-B(1-h)} + C\,e^{-D(1-h)}\Big]}.
\end{equation}
We determine $(A,B,C,D)$ by a minimax fit. The fit for $(C,D)$ was focused on high eccentricities in $\mathcal{D}_2=\{b=1,\; a\in[100,1000]\}$, while leaving the near–circular regime essentially governed by the first exponential with parameters $(A, B)$ that were also fitted. The fit was carried out with the constrain $A + C = 4.023374941 \times 10^{-4}$ which ensures eliminating the known error of $R2$ when $b/a = 0$, i.e., in extremely high eccentricities. The resulting new parameters for $\hat{P}_{\mathrm{2}}$ are:
\[
A \;\approx\; 3.37528 \times 10^{-4},
\qquad 
B \;\approx\; 10.29662,
\qquad 
C \;\approx\; 6.48093 \times 10^{-5}, 
\qquad 
D \;\approx\; 40.89043.
\]
An additional sigmoid factor $\sigma$ can be included in order to cancel the effect of the binomial of exponential terms when $a\approx b$, i.e. in low eccentricity values, in such a way that $\hat{P}_{\mathrm{2}} \rightarrow  P_{\mathrm{R2}}$, right there ($h < 0.33$) where the second Ramanujan's expression is highly accurate:
\begin{equation}
\hat{P}_{\mathrm{2}}= \frac{P_{\mathrm{R2}}}{1 - \Big[A\,e^{-B(1-h)} + C\,e^{-D(1-h)}\Big]\cdot \sigma}.
\label{eq:corr2}
\end{equation}
Here, $\sigma$ is the sigmoid function $\left ( \frac{1}{1+e^{-60(h-0.33)}} \right )$ which is completely optional and its application or absence does not alter the maximum relative error that we reported in this paper. Similarly, we have included the same sigmoidal factor in Eq. \eqref{eq:corr1} yielding $\hat{P}_{\mathrm{1}} \;=\; \frac{P_{\mathrm{R2}}}{1 - [A\,e^{-B(1-h)}] \cdot \sigma}$, so that $\hat{P}_{\mathrm{1}} \rightarrow  P_{\mathrm{R2}}$ when $a\approx b$. Thus, Eq. \eqref{eq:corr2} attains:
\begin{itemize}
\footnotesize
\item on the full range $\{b=1,\; a\in[1,\infty]\}$: maximum relative error $\approx 0.57$\,ppm,
\item on the range $\{b=1,\; a < 3.42\}$ the same error as Ramanujan's 2nd formula,
\end{itemize}
improving substantially over the one–exponential correction evaluated on the same extended range (about $32$\,ppm, worst case near $a=1000$).


\section{Experimental Results}

\subsection{Computational Setup}
All experiments use dense meshes, high–precision evaluation of $P_{\mathrm{exact}}$ via $E(\cdot)$, and reproducible code (global search bounds, seeds, and tolerances). The objective functions are evaluated pointwise on the mesh; for minimax fits, we report the supremum over the grid, and for least–squares we exploit linearity in the amplitude for fixed decay rate.
Both Eq. \eqref{eq:corr1} and Eq. \eqref{eq:corr2} are closed–form, single–line formulas that preserve numerical stability and require only elementary operations and one exponential as in Eq. \eqref{eq:corr1}, or two exponentials and a sigmoidal factor as in Eq. \eqref{eq:corr2}.

To evaluate the accuracy of the different closed-form approximations, we designed a computational framework in Python using \texttt{NumPy}, \texttt{mpmath}, and \texttt{Matplotlib}. Two dense meshes with $10^4$ uniformly distributed dots were generated to test the methods in the following way: 1) In the $h$ domain, setting $a=1000$ and varying $b\in[1, 1000]$. 2) In the $a$ domain, setting $b=1$ and varying $a\in[1.0001, 10000]$. Due to the restriction $a>b$ of Moscato's formula, we have used $1.0001$ instead of $1$.

The exact ellipse perimeter was computed using the complete elliptic integral of the second kind in Eq. (\ref{eq:exact}), whose high-precision value serves as the reference benchmark. Each approximation method (Fagnano, Euler, Ramanujan's formula 1, Ramanujan's formula 2, Cantrell, Koshy's formula 1, Koshy's formula 2, Moscato and our proposed $R2/1exp$ and $R2/2exp$) was then applied to the same mesh. For each method, we evaluated the signed and absolute value of relative error as $\;\frac{P_{\mathrm{approx}} - P_{\mathrm{exact}}}{P_{\mathrm{exact}}}\times 100\%$, and $\;\frac{|P_{\mathrm{approx}} - P_{\mathrm{exact}}|}{P_{\mathrm{exact}}}\times 100\%$ respectively. Finally, the maximum absolute value of all relative errors was extracted. 

\subsection{Relative Error Comparison (Selected Methods)}  
Figure~\ref{fig:rel_error_methods} (a) shows the signed relative error for the six most accurate closed-form approximations: Cantrell, Koshy 1, Koshy 2, Moscato, $R2/1exp$ (using $\sigma$), and $R2/2exp$ (Eq. (\ref{eq:corr2})) as a function of $a$, while $b$ remains constant and equal to $1$. From $a=1$ to $a=100$, the plot highlights that while Cantrell, Koshy 1, and Koshy 2 methods achieve errors $\approx$ $+/-0.001\%$, Moscato's and $R2/1exp$ achieve maximum relative errors of around $+/-0.00025\%$. In contrast, $R2/2exp$ remains much closer to zero and less than $+/-0.00001\%$ across the same range. From $a=100$ to $a=10000$, the error of $R2/1exp$ increases negatively from $-0.00021\%$ to $-0.004\%$ at $a=10000$, while in $R2/2exp$ method, the error remains less than $0.00001\%$ in that same interval. Figure~\ref{fig:rel_error_methods} (b) shows on a logarithmic scale the absolute relative error for the same techniques. As an illustrative reference, we locate the orbits of three known comets, all with $a < 100$.

\subsection{Summary of Maximum Errors}
Table~\ref{tab:max_errors} summarizes the maximum relative errors (in $\%$ and ppm) of ten methods computed over three different ranges: ($a=1000$ and $b\in[1,1000]$), ($b=1$ and $a\in[1.0001,100]$), and ($b=1$ and $a\in[100,1000]$). As can be seen, the oldest formulas (Fagnano and Euler) exhibit errors exceeding 10\%. The Ramanujan's approximations reduce the error drastically to the order of $0.1$\% or less. Cantrell and Koshy achieve errors in the order of $0.001$\%. On the other hand, Moscato only achieves good results, on the order of 3 ppm, when parameter $b$ is set to 1 and parameter $a$ is varied. It necessarily requires the parameter $b$ being equal to $1$  because it contains exponents where the parameter $a$ is isolated. In contrast, our proposed correction $R2/2exp$ further improves the maximum error to approximately $0.57\times 10^{-4}\%$, making it the most accurate among all closed-form, single-line known methods.

\begin{center}
\centering
\includegraphics[width=\linewidth]{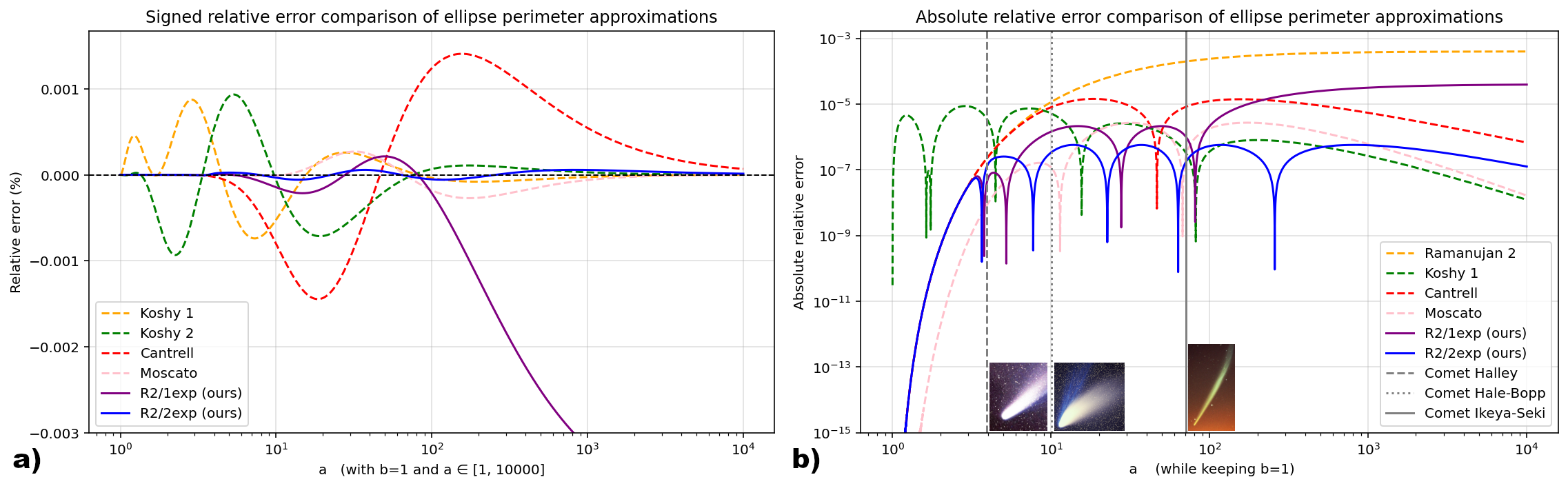}
\captionof{figure}{a) Signed relative error comparison of ellipse perimeter approximations. The graph shows that $R2/1exp$ (purple) outperforms the others methods when a<100, while $R2/2exp$ (blue) outperforms all methods across the entire range ($a<10000$). b) Absolute relative error comparison and three references of elliptical orbits of famous comets.}
\label{fig:rel_error_methods}
\end{center}

\begin{table}[h!]
\centering
\caption{Maximum relative error (\% and ppm) of ellipse perimeter approximations.}
\label{tab:max_errors}
\footnotesize 
\begin{tabular}{lccc}
\hline
\textbf{Formula} & \textbf{$a=1000$ and $b\in[1,1000]$} & \textbf{$b=1$ and $a\in[1.0001,100]$} & \textbf{$b=1$ and $a\in[100,1000]$}\\

\hline
Fagnano & 21.38195$\%$ / 213819.50 ppm & 20.69656$\%$ / 206965.61 ppm & 21.38195$\%$ / 213819.50 ppm \\
Euler & 11.0717$\%$ /110716.96 ppm & 11.04714$\%$ / 110471.35 ppm & 11.07170$\%$ / 110716.96 ppm \\
Ramanujan I & 0.4068762$\%$ / 4068.76 ppm & 0.3420281$\%$ / 3420.28 ppm & 0.4068762$\%$ / 4068.76 ppm \\
Ramanujan II & 0.03784220$\%$ / 378.42 ppm & 0.0238977$\%$ / 238.98 ppm & 0.0378422$\%$ / 378.42 ppm \\
Cantrell & 0.001446076$\%$ / 14.46 ppm & 0.001446075$\%$ / 14.46 ppm & 0.001406936$\%$ / 14.07 ppm \\
Koshy 2 & 0.0009337747$\%$ / 9.34 ppm & 0.0009337747$\%$ / 9.34 ppm & 0.0001083570$\%$ / 1.08 ppm \\
Koshy 1 & 0.0008735187$\%$ / 8.74 ppm & 0.0008735186$\%$ / 8.74 ppm & 0.00008046597$\%$ / 0.80 ppm \\
Moscato & ----- & 0.0002697398$\%$ / 2.70 ppm & 0.0002721463$\%$ / 2.72 ppm \\
\textbf{R2/1exp (proposed)} & 0.003167$\%$ / 31.67 ppm & \textbf{0.00021458$\%$} / \textbf{2.14 ppm} & 0.003167$\%$ / 31.67 ppm \\
\textbf{R2/2exp (proposed)} & \textbf{0.00005735$\%$} / \textbf{0.57 ppm} & \textbf{0.00005735$\%$} / \textbf{0.57 ppm} & \textbf{0.00005732$\%$ / 0.57 ppm} \\
\hline
\end{tabular}
\end{table}


Due to the fine-tuning of parameters $A$ and $C$, the relative error decreases from 0.57 ppm to 0 ppm as the eccentricity increases, reaching the flat ellipse limit when b/a → 0, or $e \rightarrow 1$.

\section{Discussion}

The comparative analysis of maximum relative errors provides a clear picture of the performance of the classical and modern closed-form approximations for the ellipse perimeter. The results are summarized in Table~\ref{tab:max_errors}, and the discussion below highlights both qualitative trends and quantitative comparisons. The earliest formulas, such as Fagnano and Euler, offer only rough estimates: their errors grow quickly as the ellipse becomes more eccentric, with deviations exceeding 10\%. Ramanujan's contributions mark a significant leap in accuracy: both $R1$ and $R2$ drastically reduce the error, with $R2$ already achieving errors below 0.03\%. Later refinements by Cantrell and Koshy succeeded in pushing the error below $0.01\%$ on the order of $0.001\%$, reaching an accuracy of sub-10 ppm across most of the tested range. Our exponential correction to Ramanujan's formula yields a single-line expression requiring only four exponentiations (five when the optional sigmoidal factor is included). To the best of our knowledge, this formulation achieves the lowest reported error—below 0.57 ppm across the entire range—while maintaining minimal computational overhead. Although recent approximations, such as Sykora's optimization based on Ahmadi's proposal \cite{Sykora2007Advances}, offer comparable accuracy, they involve a greater number of exponentiations, and consequently higher computational complexity. For instance, eight exponentiations in Ahmadi's proposal are used to reach $\approx 0.15$ ppm, only a quarter of ours. To contextualize the improvement, we compare the maximum relative error of each method with the proposed formula $R2/2exp$, whose maximum relative error is $0.57$ ppm, as shown in Figure \ref{fig:quantitativecomp}.
\begin{center}
\centering
\includegraphics[width=0.4\linewidth]{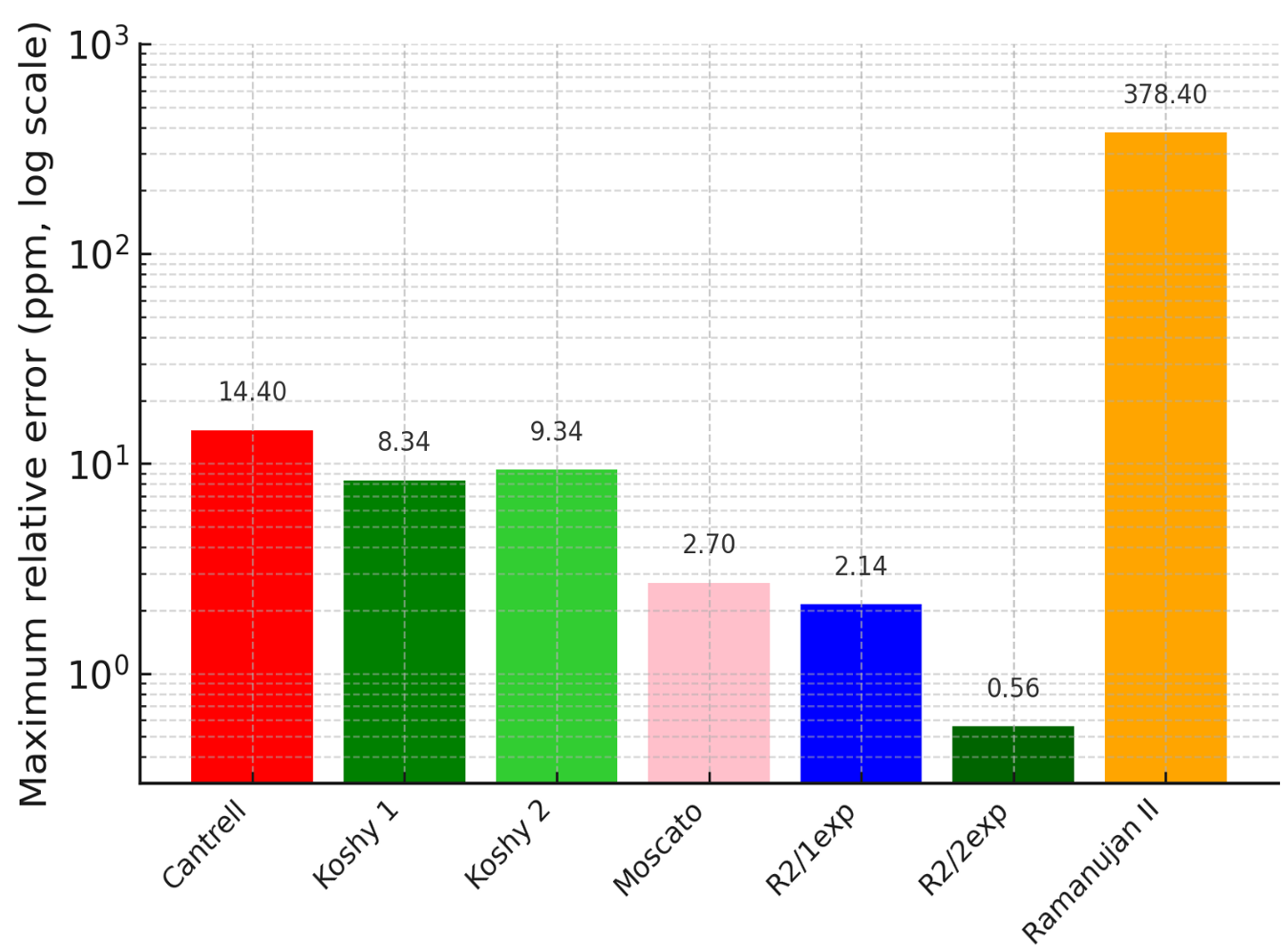}
\captionof{figure}{Closed-form ellipse perimeter approximations.}
\label{fig:quantitativecomp}
\end{center}


Finally, and as an additional piece of information we have computed Euler-Ivory expansion in Eq. (\ref{eq:Ivory}) and have seen that $148$ terms of the series are required to reach a maximum relative error of 0.57 ppm in the same interval $b=1$, and $a\in [1,1000]$ that we used for our $R2/2exp$ approximation. 

\subsection{Key insight}
The evidence demonstrates that the proposed $R2/2exp$ method consistently outperforms all other known closed-form approximations, both in magnitude and stability of the error across the full range $b/a \in [0.0, 1]$. Even the refined formulas of Cantrell, Koshy and Moscato, which are considered state of the art, still have errors between $25$, $15$, and $4.7$ times larger than ours. The resulting maximum relative error of only $0.57 \times 10^{-4}\%$ with respect to the exact perimeter makes this method the most accurate compact closed-form and single-line approximation currently available in the literature, as far as we know. One further point worth emphasizing is that although there exist certain closed-form formulas that claim maximum errors comparable or even smaller than ours, these are not single-line expressions and typically involve more fitted constants and a lot of exponents, making their expressions significantly more complex. In contrast, our proposed formula remains extraordinarily simple: it modifies Ramanujan II by dividing it by $1$ minus a binomial composed of two exponential terms, introducing just four additional parameters. This results in a formula that is short, closed-form, single-line and easily computable, achieving a maximum relative error of only about $0.57$ ppm relative to the exact perimeter of the ellipse, outperforming most other compact expressions in the trade-off between simplicity and precision. Thus, our work contributes a method with an exceptional precision-simplicity ratio: nearly the top accuracy among known formulae while preserving minimal computational and notational overhead.



\bibliographystyle{unsrt}


\bibliography{main}






\end{document}